\documentstyle{article}
\title{Between Sobolev and Poincar\'e
\thanks{Research partially supported by KBN Grant 2 P03A 043 15}}
\author{Rafa{\l} Lata{\l}a and Krzysztof Oleszkiewicz \\ {\small (Warsaw)}}
\date{}
\begin{document}
\maketitle
\newtheorem{cl}{Claim}
\newtheorem
{lem}{Lemma}
\newtheorem{cor}{Corollary}
\newtheorem{prop}{Proposition}
\newtheorem{df}{Definition}
\newtheorem{rem}{Remark}
\newtheorem{th}{Theorem}

\begin{abstract}
Let $a \in [0,1]$ and $r \in [1,2]$ satisfy relation $r=2/(2-a).$
Let $\mu(dx)=c_{r}^{n}\exp(-(|x_{1}|^{r}+|x_{2}|^{r}+\ldots +|x_{n}|^{r}))
dx_{1}dx_{2}\ldots dx_{n}$
be a probability measure on the Euclidean space $(R^{n}, \| \cdot \|).$
We prove that there exists a universal constant $C$ such that for any smooth
real function $f$ on $R^{n}$ and any $p \in [1,2)$ 
\[
E_{\mu}f^{2}-(E_{\mu}|f|^{p})^{2/p} \leq C(2-p)^{a}E_{\mu}\|\nabla f\|^{2}.
\]
We prove also that if for some probabilistic measure $\mu$
on $R^{n}$ the above inequality is satisfied for any $p \in [1,2)$ and any
smooth $f$ then for any 
$h:R^{n} \longrightarrow R$ such that $|h(x)-h(y)| \leq \|x-y\|$
there is $E_{\mu}|h| < \infty$ and
\[
\mu(h-E_{\mu}h > \sqrt{C} \cdot t) \leq e^{-Kt^{r}}
\]
for $t>1,$ where $K>0$ is some universal constant.

\end{abstract}

\bigskip
\bigskip

Let us begin with few definitions.
\begin{df}
\label{p-war}
Let $(\Omega, \mu)$ be a probability space and let $f$ be a
measurable, square integrable non-negative function on $\Omega .$
For $p \in [1,2)$ we define the $p$-variance of $f$ by
\[
Var(p)_{\mu}(f) = \int_{\Omega} f(x)^{2} \mu(dx)-
(\int_{\Omega} f(x)^{p}\mu(dx))^{2/p} =
E_{\mu}f^{2}-(E_{\mu}f^{p})^{2/p}.
\]
\end{df}

Note that $Var(1)_{\mu}(f)=D^{2}_{\mu}(f)=Var_{\mu}(f)$ coincides with
classical notion of variance, while
\[
\lim_{p \rightarrow 2^{-}} \frac{Var(p)_{\mu}(f)}{2-p} =
\frac{1}{2}(E_{\mu}f^{2} \ln (f^{2})-E_{\mu}f^{2} \cdot \ln(E_{\mu}f^{2}))
=\frac{1}{2}Ent_{\mu}(f^{2}),
\]
where $Ent_{\mu}$ denotes a classical entropy functional (see \cite{L} for 
a nice introduction to the subject).

\begin{df}
\label{nier}
Let ${\cal E}$ be a non-negative functional on some class ${\cal C}$ of
non-negative functions from $L^{2}(\Omega,\mu).$ We will say that
$f \in {\cal C}$ satisfies
\begin{itemize}
\item the Poincar\'e inequality with constant $C$\\
if $Var_{\mu}(f) \leq C \cdot {\cal E}(f),$ 
\item the logarithmic Sobolev inequality with constant $C$\\
if $Ent_{\mu}(f^{2}) \leq C \cdot {\cal E}(f),$
\item the inequality $I_{\mu}(a)$ (for $0 \leq a \leq 1$) with constant $C$\\
if $Var(p)_{\mu}(f) \leq C \cdot (2-p)^{a} \cdot {\cal E}(f)$
for all $p \in [1,2).$
\end{itemize}
\end{df}

\begin{lem}
\label{mono}
For a fixed $f \in {\cal C}$ and $p \in [1,2)$ let
\[
\varphi(p)=\frac{Var(p)_{\mu}(f)}{1/p - 1/2}.
\]
Then $\varphi$ is a non-decreasing function.
\end{lem}

{\bf Proof.}
H\"older's inequality yields that $\alpha(t)=t \ln(E_{\mu}f^{1/t})$
is a convex function for $t \in (1/2,1].$ Hence also
$\beta(t)=e^{2\alpha(t)}=(E_{\mu}f^{1/t})^{2t}$ is convex
and therefore $\frac{\beta(t)-\beta(1/2)}{t-1/2}$ is non-decreasing
on $(1/2,1].$ Observation that
\[
\varphi(p)=\frac{\beta(1/2)-\beta(1/p)}{1/p-1/2}
\]
completes the proof. $\Box$

\begin{cor}
\label{impli}
For $f \in {\cal C}$ the following implications hold true:
\begin{itemize}
\item $f$ satisfies the Poincar\'e inequality with constant $C$\\
if and only if $f$ satisfies $I_{\mu}(0)$ with constant $C,$
\item if $f$ satisfies the logarithmic Sobolev inequality with constant $C$\\
then $f$ satisfies $I_{\mu}(1)$ with constant $C,$
\item if $f$ satisfies $I_{\mu}(1)$ with constant $C$\\
then $f$ satisfies the logarithmic Sobolev inequality with constant $2C,$
\item if $f$ satisfies $I_{\mu}(a)$ with constant $C$ and $0 \leq \alpha
\leq a \leq 1$\\
then $f$ satisfies $I_{\mu}(\alpha)$ with constant $C.$
\end{itemize}
\end{cor}

{\bf Proof.}
\begin{itemize}

\item To prove the first part of Corollary \ref{impli} it suffices to note that
$p \longmapsto Var(p)_{\mu}(f)$ is a non-increasing function.

\item The second part of Corollary \ref{impli} follows easily from the fact
that 
\[\lim_{p \rightarrow 2^{-}} \frac{Var(p)_{\mu}(f)}{2-p}
= \frac{1}{2} \cdot Ent_{\mu}(f^{2}). 
\]

\item To prove the third part of Corollary \ref{impli} use Lemma \ref{mono} and
note that for $p \in [1,2)$ we have
\[
\frac{Var(p)_{\mu}(f)}{2-p} = \frac{\varphi(p)}{2p} \leq
\frac{\lim_{p \rightarrow 2^{-}} \varphi(p)}{2} =
Ent_{\mu}(f^{2}).
\]

\item The last part of statement is trivial.
$\Box$
\end{itemize}

Corollary \ref{impli} shows that inequalities $I_{\mu}(a)$ interpolate
between Poincar\'e and logarithmic Sobolev inequalities. 
Note that $I_{\mu}(a)$ for $a<0$ would be equivalent to
the Poincar\'e inequality and the only functions satisfying
$I_{\mu}(a)$ for $a>1$ would be the constant functions (because in this case
$I_{\mu}(a)$ would imply the logarithmic Sobolev inequality with constant $0$).
Therefore restriction to $a \in [0,1]$ is natural.

\begin{df}
Given probability space $(\Omega, \mu),$ a class ${\cal C} \subseteq
L^{2}_{+}(\Omega, \mu)$ and non-negative functional ${\cal E}$ on
${\cal C}$ we will say that a pair $(\mu, {\cal E})$ satisfies $I(a)$
(respectively the Poincar\'e or the logarithmic Sobolev) inequality if every
$f \in {\cal C}$ satisfies $I_{\mu}(a)$ (resp. the Poincar\'e or the
logarithmic Sobolev) inequality with constant C (for these particular $\mu$
and ${\cal E}$). For the sake of brevity we will assume that $\mu$ identifies
probability space and ${\cal E}$ carries information about ${\cal C}.$
\end{df}

An obvious modification of Corollary \ref{impli} for pairs $(\mu, {\cal E})$
follows. In some cases we can establish the precise relation between best
possible constants in $I(1)$ and logarithmic Sobolev inequalities.

Let $m: (-a,a) \longrightarrow R$ be an even, strictly postive continuous
density of some probability measure $\mu$ on $(-a,a),$ where $0 < a \leq
\infty$ and assume that $\int_{-a}^{a}x^{2}m(x)dx < \infty.$
For $f \in C_{0}^{\infty}(-a,a)$ put
\[
(Lf)(x)=xf'(x)-u(x)f''(x), 
\]
where $u(x)=\frac{\int_{x}^{a}tm(t) \ dt}{m(x)} \geq 0.$
General theory (see \cite{KLO} for detailed references and some related results)
yields that $L$ can be extended to a positive definite self-adjoint operator
(denoted by the same symbol), defined on a dense subspace $Dom(L)$ of
$L^{2}((-a,a),\mu),$ whose spectrum $\sigma(L)$ is contained in
$\{ 0 \} \cup [1,\infty).$ Moreover $P_{t}=e^{-tL}$ ($t \geq 0$) is a Markov
semigroup with invariant measure $\mu.$ Put ${\cal E}(f)=\|L^{1/2}f\|_{2}^{2}$
(we accept ${\cal E}(f)=+\infty$ for $f$ which do not belong to $Dom(L^{1/2})$)
and take ${\cal C}=L_{+}^{2}((-a,a),\mu).$

\begin{lem}
\label{L}
Under the above assumptions the following equivalence holds true:\\
$(\mu, {\cal E})$ satisifes the inequality $I(1)$ with constant $C$\\
if and only if\\
$(\mu, {\cal E})$ satisfies the logarithmic Sobolev inequality with constant
$2C.$
\end{lem}

{\bf Proof.} 
If $(\mu, {\cal E})$ satisifes the inequality $I(1)$ with constant $C$ then by 
Corollary \ref{impli} it satisfies the logarithmic Sobolev inequality
with constant $2C.$ Now let us assume that $(\mu, {\cal E})$
satisfies the logarithmic Sobolev inequality with constant $2C.$ Then for any
$f \in L^{2}((-a,a),\mu)$ we have
\[
Ent_{\mu}(f^{2})=Ent_{\mu}(|f|^{2})\leq 2C{\cal E}(|f|)\leq 2C{\cal E}(f)
\]
(the last inequality is a well known property of Dirichlet forms of Markov
semigroups - see for example Theorem 1. 3. 2 of \cite{D}).
Therefore classical hypercontractivity result \cite{G} yields
\[
\| P_{t(p)}f \|_{2} \leq \| f \|_{p},
\]
where $t(p)=\frac{C}{2}\ln(\frac{1}{p-1})$ for $p \in [1,2);$ 
if $p=1$ then we put $t(p)=\infty$ and $P_{\infty}(f)=E_{\mu}f.$
Hence
\[
Efe^{-2t(p)L}f \leq (Ef^{p})^{2/p}
\]
or equivalently
\[
Ef^{2}-(Ef^{p})^{2/p} \leq Ef(Id-e^{-2t(p)L})f
\]
for any $f \in {\cal C}.$
Now it suffices to prove that for any $\lambda \in \sigma(L)$ we have
\[
1-e^{-2t(p)\lambda} \leq (2-p)C\lambda,
\]
i.e.
\[
1-(2-p)C\lambda \leq (p-1)^{C\lambda}.
\]
For $\lambda =0$ and $p \in (1,2)$ the inequality is trivial.
It is known that if $(\mu, {\cal E})$ satisfies the logarithmic Sobolev
inequality with constant $2C$ then (under the assumptions of Lemma \ref{L})
$C \geq 1$ - to see this consider the logarithmic Sobolev inequality 
for functions of the form $f(x)=|1+\varepsilon x|$ with $\varepsilon$
tending to zero (this is a special case of more general observation which says
that, for functionals ${\cal E}$ satisfying certain natural conditions,
if $(\mu, {\cal E})$ satisfies the logarithmic Sobolev inequality with constant
$2C$ then it also satisfies the Poincar\'e inequality with
constant $C$). We can restrict our considerations to the case $\lambda \geq
1$ since $\sigma(L) \setminus \{ 0 \} \subseteq [1,\infty).$
Therefore $(p-1)^{C\lambda}$ is a convex function of $p$ and to prove
that
\[
h(p)=(p-1)^{C\lambda}+(2-p)C\lambda -1 \geq 0
\]
for $p \in [1,2)$ it suffices to check that $h(2)=h'(2)=0$
which is obvious. The case $p=1$ (omitted when $\lambda=0$ because
$(p-1)^{C\lambda}$ was not well defined) follows easily since the function
$p \longmapsto (Ef^{p})^{2/p}$ is continuous for $p \in [1,2].$
$\Box$

\begin{cor}
\label{gauss}
If $\mu$ is a ${\cal N}(0,1)$ Gaussian measure on real line,
${\cal E}(f)=E_{\mu}(f')^{2}$ and ${\cal C}$ is a class of non-negative
smooth functions then $(\mu, {\cal E})$ satisfies $I(1)$ with constant 1.
\end{cor}

{\bf Proof.}
If $\mu$ is a ${\cal N}(0,1)$ Gaussian measure and operator $L$ is defined
as before then
\[
E_{\mu}fLf=E_{\mu}(f')^{2}.
\]
The assertion follows from Lemma \ref{L} and well known fact (\cite{G}) that
Gaussian measures satisfy the logarithmic Sobolev inequality with constant 2.
$\Box$

\begin{rem}
Method used in Lemma \ref{L} seems applicable also in more general
situation (see \cite{O} for possible directions of generalization).
Let us mention just one interesting application.
If $\Omega=\{ -1,1 \},$ $\mu(\{ -1 \})=\mu(\{ 1 \})=1/2$
and ${\cal E}(f)=(\frac{f(1)-f(-1)}{2})^{2}$ then $(\mu, {\cal E})$
satisfies $I(1)$ with constant 1.
\end{rem}

\begin{rem}
 Let $\mu$ be a non-symmetric two-point distribution on $\{-1,1\}$,
 $\mu(\{ 1 \})=1-\mu(\{ -1 \})=\alpha$ with $\alpha\in (0,1/2)\cup(1/2,1)$.
 Then for any $p\in [1,2)$ and any $f:\{-1,1\}\rightarrow R_{+}$ the
 inequality
\[
E_{\mu}f^{2}-(E_{\mu}f^{p})^{2/p}\leq C_{\alpha}(p)(f(1)-f(-1))^{2}
\] 
holds with
\[
C_{\alpha}(p)=\frac{\alpha^{1-2/p}-(1-\alpha)^{1-2/p}}{\alpha^{-2/p}-(1-\alpha)^{-2/p}}
\]
and the constant cannot be improved.
\end{rem}

{\bf Proof} (sketch). To check the optimality of $C_{\alpha}(p)$ put
$f(-1)=\alpha^{2/p}$ and $f(1)=(1-\alpha)^{2/p}$. To prove the inequality
observe that for $p\in(1,2)$,
$\varphi(y)=((1+\sqrt{y})^{p}+(1-\sqrt{y})^{p})^{2/p}$ is a strictly convex
function of $y\in(0,1)$, since
\[
\varphi^{'}(y)=[(1+\sqrt{y})^{p}+(1-\sqrt{y})^{p}]^{\frac{2}{p}-1}
\frac{(1+\sqrt{y})^{p-1}-(1-\sqrt{y})^{p-1}}{\sqrt{y}}\]
\[
=\bigl(2 \sum_{k=0}^{\infty}{p\choose
2k}y^{k}\bigr)^{\frac{2}{p}-1}2\sum_{k=0}^{\infty} {{p-1}\choose {2k+1}}y^{k}
\]
is clearly increasing (note that ${p \choose 2k}$ and ${{p-1}\choose {2k+1}}$
are positive for $k=0,1,\ldots $). Hence for each $y_{0}\in (0,1)$ and $p\in
(1,2)$ there exist unique real numbers $A$ and $B$ such that
\[\varphi(y^{2})=((1+y)^{p}+(1-y)^{p})^{2/p}\geq A+By^{2} \mbox{ for all
$y\in (-1,1)$}\]
with equality holding for $|y|=y_{0}$ only. By the homogenity we may assume
that $f(-1)=(1-\alpha)^{-1/p}(1+y)$ and $f(1)=\alpha^{-1/p}(1-y)$. Putting
$y_{0}=\frac{(1-\alpha)^{1/p}-\alpha^{1/p}}{(1-\alpha)^{1/p}+\alpha^{1/p}}$,
using the above inequality after some elementary, but a little involved
computations one proves the assertion.  $\Box$

\begin{df}
\label{fi}
Let us denote by $\Phi$ the class of all continuous functions
$\varphi: [0,\infty) \longrightarrow R$ having strictly positive second
derivatve and such that $1/ \varphi ''$ is a concave function. Let us
additionally include in $\Phi$ all functions $\varphi$
of the form $\varphi(x)=ax+b,$ where $a$ and $b$ are some real constants.
\end{df}

Although it is not obvious, functions belonging to $\Phi$ form a convex
cone. There are some interesting questions connected with the class
$\Phi$ and its generalizations but we postpone them till the end of the
note.

\begin{lem}
\label{wyp}
For any $\varphi \in \Phi$ and $t \in [0,1]$ the function
$F_{t}:[0,\infty) \times [0, \infty) \longrightarrow R$ defined by
\[
F_{t}(x,y)=t\varphi(x)+(1-t)\varphi(y)-\varphi(tx+(1-t)y)
\]
is non-negative and convex.
\end{lem}

{\bf Proof.}
Non-negativity of $F_{t}$ is an easy consequence of convexity of $\varphi.$
Obviously $F_{t}$ is continuous on $[0,\infty) \times [0, \infty)$ and
twice differentiable on $(0,\infty) \times (0, \infty).$ Therefore it
suffices to prove that $Hess \, F_{t}$ (second derivative matrix) is
positive definite on $(0,\infty) \times (0, \infty).$
We skip the trivial case of $\varphi$ being an affine function.
Note that from the
positivity of $\varphi ''$ and the concavity of $1/\varphi ''$ it follows that
\[
\frac{1}{\varphi ''(tx+(1-t)y)} \geq \frac{t}{\varphi ''(x)}+
\frac{1-t}{\varphi ''(y)} \geq \frac{t}{\varphi ''(x)}.
\]
Therefore
\[
\frac{\partial^{2}F_{t}}{\partial x^{2}}(x,y)=
t\varphi ''(x)-t^{2}\varphi ''(tx+(1-t)y) \geq 0.
\]
In a similar way we prove that
$\frac{\partial^{2}F_{t}}{\partial y^{2}}(x,y) \geq 0.$
Now it is enough to prove that $\det (Hess \, F_{t}) \geq 0$
i.e. that
\[
\frac{\partial^{2}F_{t}}{\partial x^{2}}(x,y) \cdot
\frac{\partial^{2}F_{t}}{\partial y^{2}}(x,y) \geq
(\frac{\partial^{2}F_{t}}{\partial x \partial y}(x,y))^{2}
\]
which is equivalent to
\[
(t\varphi ''(x)-t^{2}\varphi ''(tx+(1-t)y))
((1-t)\varphi''(y)-(1-t)^{2}\varphi ''(tx+(1-t)y) )
\]
\[\geq
(-t(1-t)\varphi ''(tx+(1-t)y))^{2}
\]
or
\[
\varphi ''(x) \varphi ''(y) \geq
t\varphi ''(y)\varphi ''(tx+(1-t)y)+(1-t)\varphi ''(x)\varphi ''(tx+(1-t)y).
\]
After dividing by $\varphi ''(x) \varphi ''(y) \varphi ''(tx+(1-t)y)$
the last inequality follows from concavity of $1/\varphi ''$ and the proof
is complete.
$\Box$

\begin{lem}
\label{wyp2}
For a non-negative real random variable $Z$ defined on probability space
$(\Omega,\mu)$ and having finite first moment, and for $\varphi \in \Phi$
let
\[
\Psi_{\varphi}(Z)=E_{\mu}\varphi(Z)-\varphi(E_{\mu}Z).
\]
Then for any non-negative real random variables $X$ and $Y$ defined on
$(\Omega,\mu)$ and having finite first moment, and for any $t \in [0,1]$
the following inequality holds:
\[
\Psi_{\varphi}(tX+(1-t)Y) \geq t\Psi_{\varphi}(X)+(1-t)\Psi_{\varphi}(Y);
\]
in other words $\Psi_{\varphi}$ is a convex functional on the convex cone of
integrable non-negative real random variables defined on $(\Omega,\mu).$
\end{lem}

{\bf Proof.} Let us note that (under notation of Lemma \ref{wyp})
\[
\Psi_{\varphi}(tX+(1-t)Y)-t\Psi_{\varphi}(X)-(1-t)\Psi_{\varphi}(Y)=
\]
\[
(E_{\mu}\varphi(tX+(1-t)Y)-tE_{\mu}\varphi(X)-(1-t)E_{\mu}\varphi(Y))-
\]
\[
(\varphi(tE_{\mu}X+(1-t)E_{\mu}Y)-t\varphi(E_{\mu}X)-(1-t)\varphi(E_{\mu}Y))
\]
\[
=E_{\mu}F_{t}(X,Y)-F_{t}(E_{\mu}X,E_{\mu}Y)
=E_{\mu}F_{t}(X,Y)-F_{t}(E_{\mu}(X,Y)).
\]
We are to prove that it is a non-negative expression and this follows easily
from Jensen inequality. For the sake of clarity we present a detailed
argument.

Let $x_{0}=E_{\mu}X$ and $y_{0}=E_{\mu}Y.$ Lemma \ref{wyp} yields that
$F_{t}$ is convex, so that there exist constants $a,b,c \in R$ such that
\[
F_{t}(x,y) \geq ax+by+c
\]
for any $x,y \in [0,\infty)$ and
\[
F_{t}(x_{0},y_{0}) = ax_{0}+by_{0}+c.
\]
Therefore
\[
E_{\mu}F_{t}(X,Y) \geq E_{\mu}(aX+bY+c)= ax_{0}+by_{0}+c=F_{t}(x_{0},y_{0})=
F_{t}(E_{\mu}X,E_{\mu}Y)
\]
and the proof is finished.
$\Box$

\begin{lem}
\label{wyp3}
Let $(\Omega_{1},\mu_{1})$ and $(\Omega_{2},\mu_{2})$ be probability spaces
and let
$(\Omega,\mu)=(\Omega_{1} \times \Omega_{2},\mu_{1} \otimes \mu_{2})$
be their product probability space. For any non-negative random variable
$Z$ defined on $(\Omega,\mu)$ and having finite first moment 
and for any $\varphi \in \Phi$ the following inequality holds true:
\[
E_{\mu}\varphi(Z)-\varphi(E_{\mu}Z) \leq
E_{\mu}([E_{\mu_{1}}\varphi(Z)-\varphi(E_{\mu_{1}}Z)]+
[E_{\mu_{2}}\varphi(Z)-\varphi(E_{\mu_{2}}Z)]).
\]
\end{lem}

{\bf Proof.}
For $\omega_{2} \in \Omega_{2}$ let
$Z_{(\omega_{2})}$ be a non-negative random variable defined on
$(\Omega_{1},\mu_{1})$ by the formula
\[
Z_{[\omega_{2}]}(\omega_{1})=Z(\omega_{1},\omega_{2}).
\]
By Lemma \ref{wyp2} used for the probability space $(\Omega_{1},\mu_{1})$
and Jensen inequality used for the family of random variables
$(Z_{[\omega_{2}]})_{\omega_{2} \in \Omega_{2}}$
(this time we
skip the detailed argument which the reader can easily repeat after the proof
of Lemma
\ref{wyp2}) we get
\[
E_{\mu_{2}}(E_{\mu_{1}}\varphi(Z)-\varphi(E_{\mu_{1}}Z))
\geq
E_{\mu_{1}}\varphi(E_{\mu_{2}}Z)-\varphi(E_{\mu_{1}}(E_{\mu_{2}}Z))
\]
which is equivalent to the assertion of Lemma \ref{wyp3}.
$\Box$

By an easy induction argument we obtain

\begin{cor}
\label{subadd}
Let $(\Omega_{1},\mu_{1}),(\Omega_{2},\mu_{2}),\ldots,(\Omega_{n},\mu_{n})$
be probability spaces and let $(\Omega,\mu)=
(\Omega_{1}\times \Omega_{2} \times \ldots \times \Omega_{n},
\mu_{1}\otimes \mu_{2} \otimes \ldots \otimes \mu_{n})$ be their product
probability space. Let $Z$ be any integrable non-negative real random
variable defined on $(\Omega,\mu).$ 
Then for any $\varphi \in \Phi$ the following inequality holds:
\[
E_{\mu}\varphi(Z)-\varphi(E_{\mu}Z) \leq
\sum_{k=1}^{n} E_{\mu}(E_{\mu_{k}}\varphi(Z)-\varphi(E_{\mu_{k}}Z)).
\]
\end{cor}

Let us observe that the function $\varphi$ defined by $\varphi(x)=x^{2/p}$  
belongs to the class $\Phi$ if $p \in [1,2].$ Therefore by applying
Corollary \ref{subadd} to the random variable $Z=f^{p},$ where
$f \in L_{+}^{2}(\Omega,\mu),$ we obtain

\begin{cor}
\label{subadd2}
Under the notation of Corollary \ref{subadd} for any $f \in
L_{+}^{2}(\Omega,\mu)$ we have
\[
E_{\mu}f^{2}-(Ef^{p})^{2/p} \leq
\sum_{k=1}^{n} E_{\mu}(E_{\mu_{k}}f^{2}-(E_{\mu_{k}}f^{p})^{2/p}).
\]
\end{cor}

This sub-additivity property of functional $Var(p)_{\mu}$ immediately
yields the following

\begin{cor}
\label{tensor}
Assume that pairs $(\mu_{1},{\cal E}_{1}),$$(\mu_{2},{\cal E}_{2}),\ldots$
$(\mu_{n},{\cal E}_{n})$ satisfy the inequality $I(a)$ with some constant $C.$
Let $\mu=\mu_{1}\otimes \mu_{2} \otimes \ldots \otimes \mu_{n}$ and
${\cal E}(f)=E_{\mu}({\cal E}_{1}(f_{1})+{\cal E}_{2}(f_{2})+\ldots +{\cal
E}_{n}(f_{n})),$ where
\[
f_{i}(x)=f(x_{1},\ldots,x_{i-1},x,x_{i+1},\ldots,x_{n})
\]
for given $x_{1},\ldots,x_{i-1},x_{i+1},\ldots,x_{n}.$
Class ${\cal C}$ can be chosen in any way which assures that
$f \in {\cal C}$ implies $f_{i} \in {\cal C}_{i},$ for example
${\cal C}={\cal C}_{1}\otimes{\cal C}_{2}\otimes \ldots \otimes {\cal
C}_{n}.$ Then the pair $(\mu, {\cal E})$ also satisfies the inequality $I(a)$
with constant $C.$
\end{cor}

The case we will concentrate on is ${\cal E}(f)=E_{\mu}\|Ê\nabla f\|^{2}.$

\begin{prop}
\label{tensoryzacja}
Let $\mu_{1},$$\mu_{2},\ldots$$\mu_{n}$ be probability measures on $R.$
Let $C>0$ and $a \in [0,1].$ Assume that for any smooth
function $f:R\longrightarrow [0,\infty)$ the inequality
\[
E_{\mu_{i}}f^{2}-(E_{\mu_{i}}f^{p})^{2/p} \leq
C(2-p)^{a}E_{\mu_{i}}(f')^{2}
\]
holds true for $p \in [1,2)$ and $i=1,2,\ldots n.$
Then for $\mu=\mu_{1}\otimes \mu_{2} \otimes \ldots \otimes \mu_{n}$ the
inequality
\[
E_{\mu}f^{2}-(E_{\mu}f^{p})^{2/p} \leq
C(2-p)^{a}E_{\mu}\| \nabla f\|^{2},
\]
where $\| \cdot \|$ denotes standard Euclidean norm,
is satisfied for $p \in [1,2)$ and any smooth function
$f:R^{n}\longrightarrow [0,\infty).$
\end{prop}

{\bf Proof.}
Use Corollary \ref{tensor} and note that
\[
E_{\mu}\| \nabla f\|^{2}=
E_{\mu}[(\frac{\partial f}{\partial x_{1}})^{2}+\ldots+(\frac{\partial
f}{\partial x_{1}})^{2}]=
E_{\mu}[(f_{1}')^{2}+\ldots+(f_{n}')^{2}]
\]
\[
=
E_{\mu}[E_{\mu_{1}}(f_{1}')^{2}+\ldots+E_{\mu_{n}}(f_{n}')^{2}].
\, \, \, \Box
\]

Now let us demonstrate that the inequality $I(a)$ for the ${\cal
E}(f)=E_{\mu}\| \nabla f\|^{2}$ functional implies concentration of
Lipschitz functions.

\begin{th}
\label{koncentracja}
Let $\mu$ be a probability measure on $R^{n}.$ Assume that there exist
constants $C > 0$ and $a \in [0,1]$ such that the inequality
\[
E_{\mu}f^{2}-(E_{\mu} f^{p})^{2/p} \leq C(2-p)^{a}E_{\mu}\|\nabla f\|^{2}
\]
is satisfied for any smooth function $f:R^{n} \longrightarrow [0,\infty)$
and $p \in [1,2).$ Let $h:R^{n}\longrightarrow R$ be a Lipschitz function
with Lipschitz constant 1, i.e. $|h(x)-h(y)|\leq \| x-y \|$ for any
$x,y \in R^{n},$ where $\| \cdot \|$ denotes a standard Euclidean norm.
Then $E_{\mu}|h| < \infty$ and
\begin{itemize}
\item for any $t \in [0,1]$
\[
\mu(h-E_{\mu}h \geq t\sqrt{C}) \leq e^{-Kt^{2}}
\]
\item for any $t \geq 1$
\[
\mu(h-E_{\mu}h \geq t\sqrt{C}) \leq e^{-Kt^{\frac{2}{2-a}}}
\]

\end{itemize}
where $K$ is some universal constant.
\end{th}

{\bf Proof.}
Our proof will work for $K=1/3$ but we do not know optimal constants
(it is also interesting what the optimal $K$ is for given value of parameter
$a$). Note that it is essential part of the assumptions that we study the
limit behaviour when $p\rightarrow 2$. For any fixed $p\in (1,2)$ the
inequality
\[
E_{\mu}f^{2}-(E_{\mu} f^{p})^{2/p} \leq C(2-p)^{a}E_{\mu}\|\nabla f\|^{2}
\]
is weaker than the Poincar\'e inequality with constant $C(2-p)^{a}$ and
therefore it cannot imply anything stronger than the exponential
concentration.

We will follow the aproach of \cite{AS}. Assume first that
$h$ is bounded and smooth. Then $\| \nabla h \| \leq 1.$ Define
$H(\lambda)=E_{\mu}e^{\lambda h}$ for $\lambda \geq 0.$ Assumptions of Theorem
\ref{koncentracja} for $f=e^{\lambda h/2}$ give \[
H(\lambda) - H(\frac{p}{2}\lambda)^{2/p} \leq
\frac{C\lambda^{2}}{4}(2-p)^{a}E_{\mu}\| \nabla h \|^{2}e^{\lambda h} \leq
\frac{C\lambda^{2}}{4}(2-p)^{a}H(\lambda).
\]
Hence
\[
H(\lambda) \leq \frac{H(\frac{p}{2}\lambda)^{2/p}}
{1-\frac{C}{4}(2-p)^{a}\lambda^{2}}
\]
for any $p \in [1,2)$ and $0 \leq \lambda \leq
\frac{2}{\sqrt{C}}(2-p)^{-a/2}.$
Applying the same inequality for $\frac{p}{2}\lambda$ and iterating, after
$m$ steps we get
\[
H(\lambda)\leq \frac{H((\frac{p}{2})^{m}\lambda)^{(2/p)^{m}}}
{\prod_{k=0}^{m-1}(1-\frac{C\lambda^{2}}{4}(2-p)^{a}\cdot
(\frac{p}{2})^{2k})^{(2/p)^{k}} }.
\] 
Note that
\[
1-\frac{C\lambda^{2}}{4}(2-p)^{a}\cdot (\frac{p}{2})^{2k} \geq
(1-\frac{C\lambda^{2}}{4}(2-p)^{a})^{(p/2)^{2k}}
\]
since $(\frac{p}{2})^{2k} < 1.$
Hence
\[
H(\lambda)\leq H((\frac{p}{2})^{m}\lambda)^{(2/p)^{m}}
(1-\frac{C\lambda^{2}}{4}(2-p)^{a})^{-\sum_{k=0}^{m-1}(p/2)^{k}}.
\]
As $\lim_{m \rightarrow \infty} (\frac{p}{2})^{m} =0$ we get
\[
\lim_{m \rightarrow \infty} H((\frac{p}{2})^{m}\lambda)^{(2/p)^{m}}
=e^{\lambda E_{\mu}h}.
\]
Therefore
\[
E_{\mu}e^{\lambda(h-E_{\mu}h)} \leq
(1-\frac{C\lambda^{2}}{4}(2-p)^{a})^{-\frac{2}{2-p}}
\]
and
\[
\mu(h-E_{\mu}h \geq t\sqrt{C}) \leq
e^{-\lambda t\sqrt{C}} \cdot
(1-\frac{C\lambda^{2}}{4}(2-p)^{a})^{-\frac{2}{2-p}}.
\]

\begin{itemize}
\item Putting $p=1$ and $\lambda=\frac{t}{\sqrt{C}}$ we get for any
$t \in [0,2)$
\[
\mu(h-E_{\mu}h \geq t\sqrt{C}) \leq
e^{-t^{2}}\cdot (1-\frac{t^{2}}{4})^{-2}.
\]
In particular for $t \in [0,1]$ we have $1-\frac{t^{2}}{4} > e^{-t^{2}/3}$
and
\[
\mu(h-E_{\mu}h \geq t\sqrt{C}) \leq
e^{-t^{2}/3}.
\]

\item If $t\geq 1,$ let us put $p=2-t^{-\frac{2}{2-a}}$ and
$\lambda = t^{\frac{a}{2-a}}/\sqrt{C}.$ Then we arrive at
\[
\mu(h-E_{\mu}h \geq t\sqrt{C}) \leq e^{-t^{\frac{2}{2-a}}} \cdot
(1-\frac{1}{4})^{-2t^{\frac{2}{2-a}}}=(\frac{16}{9e})^{t^{\frac{2}{2-a}}}
\]
which completes the proof (if $h$ is bounded and smooth) since
$\frac{16}{9e}\leq e^{-1/3}.$
\end{itemize}

Therefore by a standard approximation argument we prove the assertion for
any bounded $h$ which satisfies assumptions of Theorem \ref{koncentracja}.
Finally for general $h$ define its bounded truncations
$(h_{N})_{N=1}^{\infty}$ putting $h_{N}(x)=h(x)$ if $|x|\leq N$ and
$h_{N}(x)=N\cdot sgn(x)$ if $|x| \geq N.$ One can easily check that
if $h$ satisfies the assumptions of Theorem \ref{koncentracja} then
$|h_{N}|$ is also a Lipschitz function with a Lipschitz constant $1$ and
therefore using Theorem \ref{koncentracja} for a bounded function $|h_{N}|$ 
we arrive at
\[
\mu(|h_{N}|-E_{\mu}|h_{N}| \geq 4\sqrt{C}) \leq
(\frac{16}{9e})^{4^{\frac{2}{2-a}}} \leq
(\frac{16}{9e})^{4} \leq \frac{1}{5}.
\]
Similarly
\[
\mu(|h_{N}|-E_{\mu}|h_{N}| \leq -4\sqrt{C}) =
\mu(-|h_{N}|-E_{\mu}(-|h_{N}|) \geq 4\sqrt{C}) \leq \frac{1}{5}.
\]
Hence
\[
\mu(|\, |h_{N}|-E_{\mu}|h_{N}| \,| \geq 4\sqrt{C}) \leq \frac{2}{5}
\]
and
\[
\mu(|\, |h|-E_{\mu}|h_{N}|\, | \geq 4\sqrt{C}) \leq \frac{2}{5}+
\mu(|h|>N).
\]
Therefore
\[
\mu(| E_{\mu}|h_{N}|-E_{\mu}|h_{M}|\, | \geq 8\sqrt{C}) \leq
\]
\[
\mu(|\, |h|-E_{\mu}|h_{N}|\, | \geq 4\sqrt{C})+
\mu(|\, |h|-E_{\mu}|h_{M}|\, | \geq 4\sqrt{C}) \leq
\]
\[
\frac{4}{5}+\mu(|h|>N)+\mu(|h|>M)
\longrightarrow \frac{4}{5} <1
\]
as $\min(N,M) \longrightarrow \infty,$
which means that the sequence $(E_{\mu}|h_{N}|)_{N=1}^{\infty}$
is bounded. As $|h_{N}|$ grows monotonically to $|h|,$ by Lebesgue
Lemma we get $E_{\mu}|h| < \infty$ and
$E_{\mu}h_{N} \longrightarrow E_{\mu}h$
as $N \longrightarrow \infty.$
Now an easy approximation argument completes the proof.
$\Box$

In order to prove that the order of concentration implied by Theorem
\ref{koncentracja} cannot be improved we will need the following
\begin{th}
\label{odwrotnenaRn}
Let $a \in [0,1]$ and $r \in [1,2]$ satisfy $r=2/(2-a).$
Put $c_{r}=\frac{1}{2\Gamma(1+1/r)}=\frac{r}{2\Gamma(1/r)}.$
Then $\mu_{r}(dx)=c_{r}^{n}\exp(-(|x_{1}|^{r}+|x_{2}|^{r}+\ldots
+|x_{n}|^{r}))dx_{1}dx_{2}\ldots dx_{n}$ is a probability measure on
$R^{n}$ and
there exists a universal constant $C>0$ (not depending on $a$ or $n$)
such that
\[
E_{\mu_{r}}f^{2}-(E_{\mu_{r}} f^{p})^{2/p} \leq
C(2-p)^{a}E_{\mu_{r}}\|\nabla f\|^{2}
\]
for any smooth non-negative function $f$ on $R^{n}$ and any $p \in [1,2).$
\end{th}

{\bf Proof.}
Proposition \ref{tensoryzacja} shows that it is enough to prove Theorem
\ref{odwrotnenaRn} in the case $n=1.$ Therefore the assertion easily
follows from the two following propositions.
$\Box$

\begin{prop}
\label{rownowaznosc}
Let $a \in [0,1]$ and $r \in [1,2]$ satisfy $r=2/(2-a).$
Put $c_{r}=\frac{1}{2\Gamma(1+1/r)},$ so that
$\mu_{r}(dx)=c_{r}\exp(-|x_{1}|^{r})dx$ is a probability measure on $R.$
Let $\lambda(dx)=\frac{1}{2}e^{-|x|}$ be a symmetric exponential
probability measure on
$R.$ Under these assumptions the following implications hold true:
\begin{itemize}
\item
If $C>0$ is a constant such that for any smooth function
$f:R\longrightarrow [0,\infty)$ and any $p \in [1,2)$
there is
\[
E_{\mu_{r}}f^{2}-(E_{\mu_{r}} f^{p})^{2/p} \leq
C(2-p)^{a}E_{\mu_{r}}(f')^{2}
\]
then for any smooth function $g:R\longrightarrow [0,\infty)$ and any $p \in
[1,2)$ there is
\[
\int_{R}g(x)^{2}\lambda(dx)-(\int_{R} g(x)^{p}\lambda(dx))^{2/p} \leq
600C(2-p)^{a}\int_{R}\max(1,|x|^{a})g'(x)^{2}\lambda(dx).
\]
\item
Conversely, if $C>0$ is such a constant that 
for any smooth function $g:R\longrightarrow [0,\infty)$ and any $p \in
[1,2)$ there is
\[
\int_{R}g(x)^{2}\lambda(dx)-(\int_{R} g(x)^{p}\lambda(dx))^{2/p} \leq
C(2-p)^{a}\int_{R}\max(1,|x|^{a})g'(x)^{2}\lambda(dx)
\]
then for any smooth function
$f:R\longrightarrow [0,\infty)$ and any $p \in [1,2)$
there is
\[
E_{\mu_{r}}f^{2}-(E_{\mu_{r}} f^{p})^{2/p} \leq
50C(2-p)^{a}E_{\mu_{r}}(f')^{2}.
\]
\end{itemize}
\end{prop}

\begin{prop}
\label{ekspo}
There exists a universal constant $C$ such that for any $a \in [0,1],$ 
any $p \in [1,2)$ and any smooth function $g:R\longrightarrow [0,\infty)$
the following inequality holds
\[
\int_{R}g(x)^{2}\lambda(dx)-(\int_{R} g(x)^{p}\lambda(dx))^{2/p} \leq
C(2-p)^{a}\int_{R}\max(1,|x|^{a})g'(x)^{2}\lambda(dx).
\]
\end{prop}

We will start with proof of Proposition \ref{rownowaznosc}. The proof of
Proposition \ref{ekspo} will be postponed to the end of the paper.

{\bf Proof of Proposition \ref{rownowaznosc}.}
Let us define the function $z_{r}:R \longrightarrow R$ by
\[
\frac{1}{2}\int_{z_{r}(x)}^{\infty}e^{-|t|}dt=
c_{r}\int_{x}^{\infty}e^{-|t|^{r}}dt,
\]
where $c_{r}=\frac{r}{2\Gamma(1/r)}=\frac{1}{2\Gamma(1+1/r)}.$
It is easy to see that $z_{r}$ is a homeomorphism of $R$ onto itself and
\[
z_{r}'(x)=2c_{r}e^{|z_{r}(x)|-|x|^{r}}.
\]
Therefore $z_{r}$ is a $C^{1}-$diffeomorphism of $R$ onto itself.
Binding $f$ and $g$ by relation $f(x)=g(z_{r}(x))$ and using standard
change of variables formula we reduce the proof of Proposition
\ref{rownowaznosc} to the following lemma.
$\Box$

\begin{lem}
\label{jakobian}
Under notation introduced above
\[
\frac{1}{50}\max(1,|x|^{a})
\leq(z_{r}'(z_{r}^{-1}(x)))^{2}
\leq 600\max(1,|x|^{a})
\]
for any $x \in R.$
\end{lem}

{\bf Proof.}
First let us note that $1/3 \leq c_{r} \leq e/2.$ Indeed,
\[
\Gamma(1/r)=\int_{0}^{\infty}x^{\frac{1}{r}-1}e^{-x}dx
 \leq
\int_{0}^{1}x^{\frac{1}{r}-1}dx+\int_{1}^{\infty}e^{-x}dx=r+1/e.
\]
Hence $c_{r}\geq \frac{r}{2r+2/e} \geq 1/3.$
On the other hand
\[
\Gamma(1/r)=\int_{0}^{\infty}x^{\frac{1}{r}-1}e^{-x}dx \geq
\frac{1}{e}\int_{0}^{1}x^{\frac{1}{r}-1}dx=r/e.
\]
Therefore $c_{r} \leq e/2.$
Let us also notice that by obvious symmetry we can consider only the case
$x>0.$
Now let us estimate from below $z_{r}^{-1}(1).$ We have
\[
\frac{e}{2}z_{r}^{-1}(1) \geq c_{r}z_{r}^{-1}(1) \geq
c_{r}\int_{0}^{z_{r}^{-1}(1)}e^{-t^{r}}dt = \frac{1}{2}\int_{0}^{1}e^{-t}dt=
\frac{1}{2}(1-1/e)
\]
and therefore $z_{r}^{-1}(1) \geq \frac{e-1}{e^{2}} \geq 1/5.$
Note that by definition of $z_{r}(x)$ for $x>0$ we have
\[
\frac{1}{2}e^{-z_{r}(x)}=c_{r}\int_{x}^{\infty}e^{-t^{r}}dt \leq
c_{r}\int_{x}^{\infty}\frac{rt^{r-1}}{rx^{r-1}}e^{-t^{r}}dt =
\frac{c_{r}e^{-x^{r}}}{rx^{r-1}}
\]
and therefore
\[
z_{r}'(x)=2c_{r}e^{z_{r}(x)-x^{r}} \geq rx^{r-1}.
\]
Hence also $z_{r}(x) \geq x^{r}$ and $z_{r}^{-1}(x) \leq x^{1/r}$
for all positive $x.$
If $x \geq 1/5$ then
\[
\int_{x}^{\infty}e^{-t^{r}}dt \geq \int_{x}^{6x}e^{-t^{r}}dt \geq
\frac{1}{r(6x)^{r-1}}\int_{x}^{6x}rt^{r-1}e^{-t^{r}}dt=
\]
\[
6^{1-r}\frac{e^{-x^{r}}-e^{-6^{r}x^{r}}}{rx^{r-1}} \geq \frac{1}{12}
\frac{e^{-x^{r}}}{rx^{r-1}},
\]
since $6^{r}x^{r} \geq x^{r}+1$ for $x \geq 1/5$ and $r \in [1,2].$
Therefore for $x \geq z_{r}^{-1}(1) \geq 1/5$ we have
\[
z_{r}'(x) \leq 12rx^{r-1} \leq 24x^{r-1}
\]
and
\[
z_{r}(x) \leq z_{r}(z_{r}^{-1}(1))+12\int_{z_{r}^{-1}(1)}^{x}rt^{r-1}dt=
1+12(x^{r}-[z_{r}^{-1}(1)]^{r}) \leq 1+12x^{r} \leq 37x^{r}.
\]
Hence $z_{r}^{-1}(x) \geq (x/37)^{1/r}$ for $x \geq z_{r}^{-1}(1).$
If $x \geq 1$ then $z_{r}^{-1}(x)\geq 1/5$ and therefore
\[
z_{r}'(z_{r}^{-1}(x)) \leq 24[z_{r}^{-1}(x)]^{r-1} \leq
24x^{\frac{r-1}{r}}=24x^{a/2}.
\]
Also if $x \geq 1$ then $z_{r}^{-1}(x) \geq z_{r}^{-1}(1)$ and
\[
z_{r}'(z_{r}^{-1}(x)) \geq r[z_{r}^{-1}(x)]^{r-1} \geq (x/37)^{\frac{r-1}{r}} \geq
37^{\frac{1}{r}-1}x^{a/2} \geq \frac{1}{7}x^{a/2}.
\]
This proves Lemma \ref{jakobian} for $|x| \geq 1.$
For any $x \geq 0$ we have
\[
z_{r}'(z_{r}^{-1}(x))=2c_{r}e^{x-z_{r}^{-1}(x)^{r}} \geq 2c_{r} \geq 2/3.
\]
We used the previously proved fact that $z_{r}^{-1}(x)\leq x^{1/r}.$
Now it remains only to establish upper estimate on $z_{r}'(z_{r}^{-1}(x))$
for $x \in [0,1].$ Note that if $x \leq z_{r}^{-1}(1)$ then
\[
c_{r}\int_{x}^{\infty}e^{-t^{r}}dt=\frac{1}{2}\int_{z_{r}(x)}^{\infty}e^{-t}dt
\geq \frac{1}{2} \int_{1}^{\infty} e^{-t}dt=\frac{1}{2e}
\]
and therefore
\[
z_{r}'(x)=\frac{2c_{r}e^{-x^{r}}}{2c_{r}\int_{x}^{\infty}e^{-t^{r}}dt} \leq
\frac{c_{r}}{c_{r}\int_{x}^{\infty}e^{-t^{r}}dt} \leq
2ec_{r} \leq e^{2} \leq 8.
\]
Hence $z_{r}'(z_{r}^{-1}(x)) \leq 8$ for any $|x| \leq 1$ and the proof is
finished.
$\Box$

\begin{lem}
\label{metr}
For $s \in (1,2]$ and $x,y \geq 0$ put
\[
\rho_{s}(x,y)=(\frac{x^{s}+y^{s}}{2}-(\frac{x+y}{2})^{s})^{1/2}.
\]
Then $\rho_{s}$ is a metric on $[0,\infty).$
\end{lem}

{\bf Proof.}
Since $k_{t}(a,b)=e^{-(a+b)t}$ is obviously positive definite integral
kernel and
$K(a,b)=s(s-1)(a+b)^{s-2}=\frac{s(s-1)}{\Gamma(2-s)}\int_{0}^{\infty}
t^{1-s}k_{t}(a,b)\, dt$  we get, by Schwartz inequality (applied to a scalar
product defined by the kernel $K(a,b)$), that for any
$y \geq x \geq 0$ and $z \geq t \geq 0$  the following inequality is true:

\begin{eqnarray*}
\int_{x/2}^{y/2} \int_{t/2}^{z/2} && \!\!\!\!\!\!\!\!\! K(a,b) \, da \, db
\\ &\leq &
(\int_{x/2}^{y/2} \int_{x/2}^{y/2} K(a,b) \, da \, db)^{1/2}
(\int_{t/2}^{z/2} \int_{t/2}^{z/2} K(a,b) \, da \, db)^{1/2}. 
\end{eqnarray*}

Now, as

\[
K(a,b)=\frac{\partial^{2}}{\partial a \, \partial b} (a+b)^{s},
\]
we get by integration by parts

\[
(\frac{y+z}{2})^{s}+(\frac{x+t}{2})^{s}
-(\frac{x+z}{2})^{s}-(\frac{y+t}{2})^{s}
\leq
\]

\[(x^{s}+y^{s} -2(\frac{x+y}{2})^{s})^{1/2}
  (z^{s}+t^{s} -2(\frac{z+t}{2})^{s})^{1/2}\]
Putting $t=y$ we arrive at

\[
(\frac{x+y}{2})^{s}+(\frac{y+z}{2})^{s}
-(\frac{x+z}{2})^{s}-y^{s}
\leq
2\rho_{s}(x,y)\rho_{s}(y,z)
\]
which is equivalent to

\[
\rho_{s}(x,z)^{2}-\rho_{s}(x,y)^{2}-\rho_{s}(y,z)^{2}
\leq
2\rho_{s}(x,y)\rho_{s}(y,z).
\]
Hence $\rho_{s}(x,z)\leq \rho_{s}(x,y)+\rho_{s}(y,z).$ For $x \leq y \leq z$
we have also easily $\rho_{s}(x,z) \geq \rho_{s}(x,y)$ and
$\rho_{s}(x,z) \geq \rho_{s}(y,z),$ so that
$\rho_{s}(x,y) \leq \rho_{s}(x,z)+\rho_{s}(z,y)$ and
$\rho_{s}(y,z) \leq \rho_{s}(y,x)+\rho_{s}(x,z).$ 
This completes the proof of triangle inequality for $s<2.$ Other metric
properties of $\rho_{s}$ as well as the case $s=2$ are trivial. $\Box$

\begin{rem}
In a similar way one can prove that
$\rho_{s}(x,y)=|\frac{x^{s}+y^{s}}{2}-(\frac{x+y}{2})^{s}|^{1/2}$
is a metric on $(0, \infty)$ for $s \in (-\infty,0) \cup (0,1).$ It was pointed
out to the authors by B. Maurey that Lemma \ref{metr} follows also from
isometrical
immersion of $([0,\infty),\rho_{s})$ into
$L^{2}([0,\infty),\kappa_{s}^{-1}t^{-s-1}dt),$ where $x\in [0,\infty)$ is sent
to the function $e^{-xt}-1$ and
$\kappa_{s}=2^{s+1}\int_{0}^{\infty}(e^{-u}-1+u)u^{-s-1}du.$ 
\end{rem}

\begin{lem}
\label{cl}
Let $s\in[1,2]$, $t\in [0,1]$ and $c,d,x$ be nonnegative numbers.
The following inequality holds
\[(1-t)c^{s}+td^{s}-((1-t)c+td)^{s}\leq \] 
\begin{equation}
\label{claim}  
  K[(1-t)c^{s}+td^{s}+x^{s}-((1-t)c+tx)^{s}-(td+(1-t)x)^{s}].
\end{equation}
under anyone of the following additional assumptions
\begin{itemize}
  \item $x$ lies outside the open interval $(c,d)$ and $K=1$
  \item $t=\frac{1}{2}$ and $K=2$
  \item $t\leq \frac{1}{2}$, $c\geq d$ and $K=12$
\end{itemize}
\end{lem}

{\bf Proof.}
Let us remind  that
\[F_{t}(x,y)=tx^{s}+(1-t)y^{s}-(tx+(1-t)y)^{s}\] is a convex function
on $[0,\infty) \times [0,\infty).$ Note that the inequality of Lemma \ref{cl}
is equivalent to \[F_{t}(d,c) \leq K[F_{t}(d,x)+F_{t}(x,c)].\] 

\begin{itemize}

\item
As
\[\frac{\partial}{\partial x}[F_{t}(d,x)+F_{t}(x,c)]\]
\[=s[(1-t)(x^{s-1}-(td+(1-t)x)^{s-1})+t(x^{s-1}-(tx+(1-t)c)^{s-1})],\]
we see that the right-hand side of the inequality as a function
of $x$ is increasing on $(\max(c,d), \infty)$ and decreasing on
$[0, \min(c,d)).$ For $x=\max(c,d)$ and $x=\min(c,d)$ the inequality is
trivially satisfied with $K=1.$ This completes the case of $x$ which does
not lie between $c$ and $d.$

\item
The second part of Lemma \ref{cl} follows easily by Lemma \ref{metr}, as
\[
F_{1/2}(d,c)=\rho_{s}(d,c)^{2} \leq (\rho_{s}(d,x)+\rho_{s}(x,c))^{2} \leq
\]
\[
2[\rho_{s}(d,x)^{2}+\rho_{s}(x,c)^{2}]=2[F_{1/2}(d,x)+F_{1/2}(x,c)].
\]

\item
To prove the last part of the statement we will use convexity of $F_{t}.$
Since $F_{t}(d,x)+F_{t}(x,c) \geq F_{t}(\frac{d+x}{2},\frac{x+c}{2}),$
it suffices to prove that $F_{t}(d,c) \leq
12 F_{t}(\frac{d+x}{2},\frac{x+c}{2}).$ Thanks to the first part of Lemma 
\ref{cl} we
can restrict our considerations to the case $x \in [d,c].$ Note that
\[\frac{\partial}{\partial x} F_{t}(\frac{d+x}{2},\frac{x+c}{2})\] 
\[=\frac{s}{2}[t(\frac{d+x}{2})^{s-1}+(1-t)(\frac{x+c}{2})^{s-1}
-(t(\frac{d+x}{2})+(1-t)(\frac{x+c}{2}))^{s-1}] \leq 0,
\]
since the function $\varphi (u)=u^{s-1}$ is concave. Therefore it is enough
to prove that
\[
F_{t}(d,c) \leq 12 F_{t}(\frac{d+c}{2},c).
\]
Using the homogenity of the above formula we can reduce our task to proving
that \[
F_{t}(1-u,1) \leq 12 F_{t}(1-u/2,1)
\]
for any $u \in [0,1]$ and $t \in [0,1/2].$

Using the Taylor expansion we get

\[
F_{t}(1-u,1)=t(1-u)^{s}+1-t-(1-tu)^{s}=
\]
\[
s(s-1)u^{2}t(1-t) \cdot \bigl[\frac{1}{2} +
\sum_{k=1}^{\infty} \frac{u^{k}}{(k+1)(k+2)}
\sum_{m=0}^{k} t^{m} \cdot \prod_{l=1}^{k}(1-\frac{s-1}{l})\bigr].
\]
Therefore
\[
F_{t}(1-u/2,1) \geq \frac{1}{2}s(s-1)(u/2)^{2}t(1-t)
\]
and
\[
F_{t}(1-u,1) \leq s(s-1)u^{2}t(1-t) \cdot \bigl[\frac{1}{2} + 
2\sum_{k=1}^{\infty} \frac{1}{(k+1)(k+2)}\bigr]\]
\[=\frac{3}{2}s(s-1)u^{2}t(1-t)\]
because $\sum_{m=0}^{\infty} t^{m} \leq 2$.
Hence
\[F_{t}(1-u,1) \leq 12 F_{t}(1-u/2,1)\]
which completes the proof. $\Box$
\end{itemize}

\begin{lem}
\label{Schw}
  Let $a\in [0,1]$, $0\leq x_{1}<x_{2}$ and $g$ be a smooth function on
  $[x_{1},x_{2}]$ such that $g(x_{1})=y_{1}, g(x_{2})=y_{2}$. Then
  \begin{equation}
  \label{est1}
   \int_{x_{1}}^{x_{2}}\max(1,x^{a})g'(x)^{2}d\lambda(x)\geq
   \frac{(y_{2}-y_{1})^{2}}{4(e^{x_{2}}-e^{x_{1}})}\max(1,x_{2}^{a}).
  \end{equation}
 
\end{lem}

{\bf Proof.} By the Schwartz inequality
\[|y_{2}-y_{1}|\leq \int_{x_{1}}^{x_{2}}|g'(x)|dx\]
\[\leq(\int_{x_{1}}^{x_{2}}\max(1,x^{a})g'(x)^{2}d\lambda(x))^{1/2}
  (2\int_{x_{1}}^{x_{2}}\min(1,x^{-a})e^{x}dx)^{1/2}.\]
Therefore to show (\ref{est1}) it is enough to prove that
\[f_{1}(x_{2})=\int_{x_{1}}^{x_{2}}\min(1,x^{-a})e^{x}dx\leq
  2\min(1,x_{2}^{-a})(e^{x_{2}}-e^{x_{1}})=f_{2}(x_{2}).\] 
For $x_{2}\leq 2$ this is obvious
because for $0<x<x_{2}\leq 2$ we have
$\min(1,x^{-a})\leq 1 \leq 2\min(1,x_{2}^{-a}),$
and for $x\geq 2$ we have
\[f_{2}'(x)=2x^{-a}(e^{x}-ax^{-1}(e^{x}-e^{x_{1}}))\geq
x^{-a}e^{x}
=f_{1}'(x).
  \Box \]

\begin{lem}
\label{Schwinv}
 Let $0\leq y_{1}<y_{2}$, $0\leq x_{1}<x_{2}$ and $g$ is defined on
 $(-\infty,x_{2})$ by the formula 
  \[g(x)=\left\{\begin{array}{lc}
               y_{1}&\mbox{for $x\leq x_{1}$}\\
               y_{1}+(e^{x}-e^{x_{1}})\frac{y_{2}-y_{1}}{e^{x_{2}}-e^{x_{1}}}&
               \mbox{for $x\in (x_{1},x_{2}]$}
               \end{array}\right. .\]
  Then
  \begin{equation}
  \label{est2}
   \int_{-\infty}^{x_{2}}g'(x)^{2}d\lambda(x)=
   \frac{(y_{2}-y_{1})^{2}}{2(e^{x_{2}}-e^{x_{1}})}.
  \end{equation}
  and for all $p\geq 1$
  \begin{equation}
  \label{est3}
   \int_{-\infty}^{x_{2}}g(x)^{p}d\lambda(x)\leq \lambda(-\infty,x_{2})
   [(1-\frac{x_{2}}{2}e^{-x_{2}})y_{1}^{p}+\frac{x_{2}}{2}e^{-x_{2}}y_{2}^{p}].
  \end{equation}
\end{lem}

{\bf Proof.} Equation (\ref{est2}) follows by direct calculations.
It is easy to see that $g(x)$ is maximal (for fixed values of $x_{2}, y_{1}$ 
and $y_{2}$) when $x_{1}=0$, so to prove (\ref{est3}) we may and will assume
that this is the case. To easy the notation we will denote $x_{2}$ by $x$.
First we will consider $p=1$. After some standard calculations (\ref{est3})
is equivalent in this case to
\[\frac{e^{x}(x-1+e^{-x})}{(2e^{x}-1)(e^{x}-1)}\leq \frac{1}{2}xe^{-x}
  \mbox{ for all $x>0$},\]
that is 
\[2+3x\leq xe^{-x}+2e^{x} \mbox{ for all $x>0$},\]
which immeditely follows from well known estimates $e^{-x}\geq 1-x$ and
$e^{x}\geq 1+x+x^{2}/2$. 

Now, for arbitrary $p\geq 1$ notice that
$g(x)=(1-\theta(x))y_{1}+\theta(x)y_{2}$  with
$0\leq \theta(x)\leq 1.$ Therefore we have by the convexity of $x^{p}$
\[ \int_{-\infty}^{x_{2}}g(x)^{p}d\lambda(x)\leq 
\int_{-\infty}^{x_{2}}((1-\theta(x))y_{1}^{p}+
\theta(x)y_{2}^{p})d\lambda(x)\leq
\]
\[\lambda(-\infty,x_{2})
   [(1-\frac{x_{2}}{2}e^{-x_{2}})y_{1}^{p}+
   \frac{x_{2}}{2}e^{-x_{2}}y_{2}^{p}],\]
where the last inequality follows by the previously established case
$p=1$.
$\Box$

\begin{lem}
\label{add}
 Suppose that $s\in (1,2]$, $t\in (0,1)$, 
 $u=\frac{s}{4(s-1)}e^{-s/2(s-1)}$ and positive numbers 
 $a,b,c,d,\tilde{a},\tilde{c},x$ satisfy the following 
 conditions
 \[c<x<d, c^{s}\leq a, d^{s}\leq b, \tilde{c}^{s}\leq \tilde{a},
   \tilde{c}\leq (1-u)c+ux.\]
 Then 
 \[(1-t)a+tb-((1-t)c+td)^{s}\leq\]
 \begin{equation}
 \label{cl2} 
  8[(1-t)\tilde{a}+tb-((1-t)\tilde{c}+td)^{s}+(1-t)a+tx^{s}-((1-t)c+tx)^{s}].
 \end{equation}  
\end{lem}

{\bf Proof.} Without loss of generality we may assume that
$a=c^{s},b=d^{s}, \tilde{a}=\tilde{c}^{s}$. Since the function
$y\rightarrow (1-t)y^{s}-((1-t)y+td)^{s}$ is nonincreasing on $[0,d]$, it is
enough to show that
\[(1-t)c^{s}+td^{s}-((1-t)c+td)^{s}\leq\]
\[3[(1-t)((1-u)c+ud)^{s}+td^{s}-((1-t)(1-u)c+(t+(1-t)u)d)^{s}].\]
By the homogenity we may and will assume that $d=1$. 
We are then to show that
\begin{equation}
\label{cl3}
  f((1-c))\leq 8f((1-u)(1-c)),
\end{equation}  
where
\[
f(x)=(1-t)(1-x)^{s}+t-(1-(1-t)x)^{s}=\sum_{i=2}^{\infty}(-1)^{i}
{s\choose i}(1-t)(1-(1-t)^{i-1})x^{k}.
\]
We use the following simple observation: if $a_{i},b_{i}$ are two summable
sequences of positive numbers such that for any $i>j$, $a_{i}/a_{j}\geq
b_{i}/b_{j}$ then for any nondecreasing nonnegative sequence $c_{i}$
\[\frac{\sum a_{i}c_{i}}{\sum a_{i}}\geq \frac{\sum b_{i}c_{i}}{\sum b_{i}}.\]
We apply the above to the sequences $a_{i}=(-1)^{i}{s\choose i}
(1-t)(1-(1-t)^{i-1})x^{i}$, $b_{i}=(i-1)(-1)^{i}{s\choose i}$ and
$c_{i}=(1-u)^{i}$, $i=2,3,\ldots $ and notice that
\[
h(y):=\sum_{i=2}^{\infty}b_{i}y^{i}=1-(1-y)^{s-1}(1+(s-1)y) \mbox{ for $y\in
[0,1]$}
\]  
Therefore we get
\[f((1-u)x)\geq \frac{h(1-u)}{h(1)}=\bigl(1-u^{s-1}(1+(s-1)(1-u))\bigr)f(x)\]
Inequality (\ref{cl3}) follows if we notice that 
\[u^{s-1}(1+(s-1)(1-u))\leq
su^{s-1}=\frac{s^{s}}{4^{s-1}}e^{-s/2}(\frac{1}{s-1})^{s-1}\leq
1e^{-1/2}e^{1/e}\leq \frac{7}{8} \Box\]
  
\begin{prop}
\label{wynikanie}
Suppose that for all $p\in [1,2)$ and all nonnegative smooth functions $g$ we 
have
\begin{equation}
\label{ext}
  \int_{R}g^{2}d\lambda-(\int_{R}g^{p}d\lambda)^{2/p}\leq
   K_{1}(2-p)^{i}\int_{R}(g'(x))^{2}\max(1,|x|^{i})d\lambda(x) 
   \mbox{ for } i=0,1,
\end{equation}   
where $K_{1}$ is a universal constant. Then for all $p$ and $g$
as above we have
\[\int_{R}g^{2}d\lambda-(\int_{R}g^{p}d\lambda)^{2/p}\leq\]
\begin{equation}
\label{alla}
   K_{2}(2-p)^{a}\int_{R}(g'(x))^{2}\max(1,|x|^{a})d\lambda(x) 
   \mbox{ for } a\in(0,1),
\end{equation}   
where $K_{2}\leq 32K_{1}$ is some universal constant.
\end{prop}

{\bf Proof.} An easy approximation argument shows that (\ref{ext}) holds 
for any continuous function $g$, continuously differentiable everywhere
except possibly finitely many points.

First we assume that $g$ is constant on $R^{-}$ or $R^{+}$, without loss
of generality say it is $R^{-},$ and we show that (\ref{alla}) holds 
with $K_{2}=16K_{1}.$
Let us fix $p\in [1,2)$ and define 
\[x_{p}=(2-p)^{-1}, y=g(x_{p}), t=\lambda(x_{p},\infty), s=\frac{2}{p},\]
\[a=\frac{1}{1-t}\int_{-\infty}^{x_{p}}g^{2}d\lambda,
  b=\frac{1}{t}\int_{x_{p}}^{\infty}g^{2}d\lambda\]
\[c=\frac{1}{1-t}\int_{-\infty}^{x_{p}}g^{p}d\lambda \mbox{ and }
  d=\frac{1}{t}\int_{x_{p}}^{\infty}g^{p}d\lambda.\]
Notice that by H\"older's inequality we have
\begin{equation}
\label{hold}
 a\geq c^{s} \mbox{ and } b\geq d^{s}.
\end{equation}  
  
We will consider two cases

{\bf Case 1}. $y^{p}$ lies outside $(c,d)$ or $c>d$.

We first apply inequality (\ref{ext}) for $i=1$ and a function
$gI_{(-\infty,x_{p})}+yI_{[x_{p},\infty)}$ to get
\[(1-t)a+ty^{2}-((1-t)c+ty^{p})^{s}\leq 
  K_{1}(2-p)\int_{0}^{x_{p}}(g'(x))^{2}\max(1,|x|)d\lambda(x)\leq\]
\[K_{1}(2-p)^{a}\int_{0}^{x_{p}}(g'(x))^{2}\max(1,|x|^{a})d\lambda(x).\]
In a similar way using the case of $i=0$ for the function 
$yI_{(-\infty,x_{p})}+gI_{[x_{p},\infty)}$ we get
\[tb+(1-t)y^{2}-(td+(1-t)y^{p})^{s}\leq 
  K_{1}\int_{x_{p}}^{\infty}(g'(x))^{2}d\lambda(x)\leq\]
\[K_{1}(2-p)^{a}\int_{x_{p}}^{\infty}(g'(x))^{2}\max(1,|x|^{a})d\lambda(x).\]
Notice also that
\[\int_{R}g^{2}d\lambda-(\int_{R}g^{p}d\lambda)^{2/p}=
  (1-t)a+tb-((1-t)c+td)^{s}\leq\]
\[12\bigl[(1-t)a+ty^{2}-((1-t)c+ty^{p})^{s}+tb+(1-t)y^{2}-(td+(1-t)y^{p})^{s}
  \bigr]\leq\]
\[12K_{1}(2-p)^{a}\int_{R}(g'(x))^{2}\max(1,|x|^{a})d\lambda(x).\]  
The middle inequality follows by Lemma \ref{cl} with $x=y^{p}$ together
with estimates (\ref{hold}).

{\bf Case 2}. $c<y^{p}<d$, we can then find $0< x_{0}<x_{p}$ such that
$g(x_{0})=c^{1/p}$. Define new function $f$ by the formula
\[f(x)=\left\{\begin{array}{lc}
               g(x)&\mbox{for $x>x_{p}$}\\
 	       c^{1/p}+\frac{y-c^{1/p}}{e^{x_{p}}-e^{x_{0}}}(e^{x}-e^{x_{0}})
 	        &\mbox{for $x\in [x_{0},x_{p}]$}\\
 	       c^{1/p}&\mbox{for $x<x_{0}$}.
                       \end{array}\right.\]
Let
\[\tilde{a}=\frac{1}{1-t}\int_{-\infty}^{x_{p}}f^{2}d\lambda \mbox{ and }
  \tilde{c}=\frac{1}{1-t}\int_{-\infty}^{x_{p}}f^{p}d\lambda.\]                      
By Lemma \ref{Schw} and \ref{Schwinv} we have
\[\int_{R} f'(x)^{2}d\lambda(x)\leq 2(2-p)^{a}\int_{R}\max(1,|x|^{a})g'(x)^{2}
  d\lambda(x).\] 
Therefore by (\ref{ext}) with $i=0$, used for the function $f$, we have
\[(1-t)\tilde{a}+tb-((1-t)\tilde{c})+td)^{s}\leq
  2K_{1}(2-p)^{a}\int\max(1,|x|^{a})g'(x)^{2}d\lambda(x).\]
We conclude as in the previous case using Lemmas \ref{Schwinv} and \ref{add}
instead of Lemma \ref{cl}.

Finally suppose that $g$ is arbitrary. A similar argument as in case 1 (but now
with $x_{p}=0$ and $t=1/2$) together with the already proved case of $g$
constant on $R_{-}$ or $R_{+}$ proves the assertion in this case. $\Box$

\bigskip

{\bf Proof of Proposition \ref{ekspo}.}
We need only to prove that assumptions of Proposition \ref{wynikanie} are
satisfied. But in view of Proposition \ref{rownowaznosc} they are equivalent
to the Poincar\'e inequality for symmetric exponential probability measure
($i=0$) and the logarithmic Sobolev inequality for the centered
${\cal N}(0,\sqrt{2}/2)$ Gaussian measure
($i=1$) which are well known to hold with some universal constants.
This completes the proof.
$\Box$

\bigskip
In the end of the paper we would like to come back to the class $\Phi$
introduced in Definition \ref{fi}. It is easy to check that if Lemma \ref{wyp3}
holds for some function $\varphi \in C^{2}((0,\infty))\cap C([0,\infty))$
for any $(\Omega_{1},\mu_{1}),
(\Omega_{2},\mu_{2})$ and any $Z$ then $\varphi \in \Phi.$ Indeed, it is even
true if we restrict our consideration to $(\Omega_{1},\mu_{1})$ and
$(\Omega_{2},\mu_{2})$ being two-point probability spaces whose atoms have
$1/2$ measures. This gives a natural characterization of the class $\Phi.$

One can try to generalize the definition of $\Phi.$ Let $U$ be an open, convex
subset of $R^{d}.$ 
We will say that a continuous function $f:U \longrightarrow R$
belongs to the class $C_{n}(U)$ if for any probability spaces
$(\Omega_{1},\mu_{1}), \ldots ,(\Omega_{n},\mu_{n})$ and any integrable random
variable $Z$ with values in $U,$ defined on $(\Omega,\mu)=
(\Omega_{1} \times \ldots \times \Omega_{n}, \mu_{1} \otimes \ldots
\otimes \mu_{n})$ the following inequality is satisfied:
\[
\sum_{K \subseteq \{ 1,2,\ldots ,n \} } (-1)^{|K|}E_{K^{c}}f(E_{K}Z) \geq 0,
\]
where $E_{K}$ denotes expectation with respect to $\mu_{k}$ for all $k \in
K.$
One can easily see that $C_{1}(U)$ is just a set of all convex functions
on $U,$ while $C_{2}((0,\infty))$ is closely related to the class $\Phi.$
In fact $f \in C_{2}((0,\infty))$ if and only if it is an affine
function or it has a continuous strictly positive second derivative such
that $1/f''$ is a concave function. One can prove that always
$C_{n+1}(U) \subseteq C_{n}(U)$ and therefore it is natural to define
$C_{\infty}(U)$ as an intersection of all $C_{n}(U).$ Then it appears that
$f \in C_{\infty}(U)$ if and only if $f$ is given by the formula
$f(x)=Q(x,x)+x^{*}(x)+y,$ where $Q$ is a non-negative definite symmetric
quadratic form, $x^{*}$ is a linear functional on $R^{d}$ and $y$ is a
constant. The above inclusions do not need to be strict. For example it is
easy to see that $C_{2}(R)$=$C_{\infty}(R).$ It would be interesting to know
some nice characterization of $C_{2}(U)$ for general $U$ and
$C_{n}((0,\infty))$ for $n>2.$ It is not clear what applications of $C_{n}$
for $n>2$ could be found but it is easy to see that this class has some
tensorization property. By now, we do not know even the answer to the
following question: For which $p \in [1,2]$ does $f(x)=x^{p}$
belong to $C_{n}((0,\infty))?$ We can only give some estimates.

These problems will be discussed in a separate paper.

\begin{rem}
Recently some new results were announced to the authors by F. Barthe
(private communication) - he proved (using Theorem 2 above) that if a
log-concave probability measure $\mu$ on the Euclidean space $(R^{n},\| \cdot
\|)$ satisfies inequality $\mu(\{ x\in R^{n}\, ; ||x||>t \}) \leq
ce^{-(t/c)^{r}}$ for some constants $c>0, r\in [1,2]$
and any $t>0$ then it satisfies also inequality
\[
E_{\mu}f^{2}-(E_{\mu}f^{p})^{2/p} \leq C(c,n,r)(2-p)^{a}E_{\mu}\| \nabla f\|^{2}
\]
for any non-negative smooth function $f$ on $R^{n}$ and $p\in[1,2),$
where $C(c,n,r)$ is some positive constant depending on $c,n$ and $r$ only 
and $a=2-2/r$.
\end{rem}

\bigskip

{\bf Acknowledgements.} The article was inspired by the questions of Prof.
Stanis{\l}aw Kwapie\'n and Prof. Gideon Schechtman. This work was done while
the first named author was     visiting Southeast Applied Analysis Center at
School of Mathematics,     Georgia Institute of Technology and was partially
supported by NSF     Grant  DMS 96-32032. The research of    the second named
author was performed at      the Weizmann Institute of Science
   in Rehovot, Israel and Equipe d'Analyse, Universit\'e Paris VI.

\sc \noindent
Institute of Mathematics\\
Warsaw University\\
Banacha 2\\
02-097 Warszawa\\
Poland\\
E-mail: rlatala@mimuw.edu.pl, koles@mimuw.edu.pl\\
\end{document}